\documentclass[letterpaper,oneside,10pt]{article}
\usepackage[USenglish]{babel}
\usepackage[T1]{fontenc}
\usepackage[ansinew]{inputenc}

\usepackage{lmodern} %Type1-font for non-english texts and characters

\usepackage{graphicx} %%For loading graphic files
\usepackage{float}
\usepackage{caption}
\usepackage{amsmath}
\usepackage{amsthm}
\usepackage{amsfonts}
\usepackage{amssymb}
\usepackage{mathtools}
\usepackage[hyphens]{url}

\let\Right\right 
\let\Left\left 
\makeatletter 
\def\right#1{\Right#1\@ifnextchar){\!\right}{}} 
\def\left#1{\Left#1\@ifnextchar({\!\left}{}} 
\makeatother

\begin{document}

\pagestyle{empty} %No headings for the first pages.

%% Title Page %%%%%%%%%%%%%%%%%%%%%%%%%%%%%%%%%%%%%%%%%%%%%%%
%% ==> Write your text here or include other files.

%% The simple version:
\title{Calculations relating to some \\ special Harmonic numbers}
\author{John Blythe Dobson (j.dobson@uwinnipeg.ca)}

%\date{} %%If commented, the current date is used.
\maketitle

%% Inhaltsverzeichnis %%%%%%%%%%%%%%%%%%%%%%%%%%%%%%%%%%%%%%%
% \tableofcontents %Table of contents
% \cleardoublepage %The first chapter should start on an odd page.

\pagestyle{plain} %Now display headings: headings / fancy / ...

%%%%%%%%%%%%%%%%%%%%%%%%%%%%%%%%%%%%%%%%%%%%%%%%%%%%%%%%%%%%%
%% ==> Some hints are following:

\begin{abstract}
\noindent
We report on the results of a computer search for primes $p$ which divide an Harmonic number $H_{\lfloor p/N \rfloor}$ with small $N > 1$.

\end{abstract}

\section{Introduction}

\noindent
Before its complete proof by Andrew Wiles, a major result for the first case of Fermat's last theorem (FLT), that is, the assertion of the impossibility of $x^p + y^p = z^p$ in integers $x, y, z$, none of which is divisible by a prime $p > 2$, was the 1909 theorem of Wieferich \cite{Wieferich} that the exponent $p$ must satisfy the congruence

\begin{displaymath}
q_p(2) := \frac{2^{p-1} - 1}{p} \equiv 0 \pmod{p},
\end{displaymath}

\noindent
where $q_p(2)$ is known as the Fermat quotient of $p$, base 2. This celebrated result was to be generalized in many directions. One of these was the extension of the congruence to bases other than 2, the first such step being the proof of an analogous theorem for the base 3 by Mirimanoff \cite{Mirimanoff} in 1910. Another grew from the recognition that some of these criteria could be framed in terms of certain special Harmonic numbers of the form

\begin{equation} \label{eq:Harmonic}
H_{\lfloor p/N \rfloor} := \sum_{j=1}^{\lfloor p/N \rfloor} \frac{1}{j}
\end{equation}

\noindent
for $N > 1$, and with $\lfloor \cdot \rfloor$ denoting the greatest-integer function. There are nice historical overviews of these developments in \cite{Ribenboim1978} and \cite{Ribenboim1979} (pp.\ 155--59), which go into greater detail than we attempt here.

\section{Congruences for $H_{\lfloor p/N \rfloor}$ involving only Fermat quotients and low-order linear recurrent sequences}

\noindent
At the time when Wieferich's and Mirimanoff's results appeared, it was already known that three of these Harmonic numbers had close connections to the Fermat quotient, and satisfied the following congruences (all modulo $p$):

\begin{align} \label{eq:P2}
H_{\lfloor p/2 \rfloor} & \equiv -2 \cdot q_p(2) \\
H_{\lfloor p/3 \rfloor} & \equiv -\frac{3}{2} \cdot q_p(3) \\
H_{\lfloor p/4 \rfloor} & \equiv -3 \cdot q_p(2).
\end{align}

\noindent
All these results are due to Glaisher; those for $N$ = 2 and 4 will be found in (\cite{Glaisher}, pp.\ 21-22, 23), and that for $N$ = 3 in (\cite{Glaisher3}, p.\ 50). Although it was not the next Harmonic number criterion published, it will be convenient to dispense next with the case $N$ = 6. The apparatus needed to evaluate this Harmonic number appears in a 1905 paper of Lerch (\cite{Lerch}, p.\ 476, equations 14 and 15), but the implications of Lerch's result were long overlooked, and only realized in 1938 by Emma Lehmer (\cite{Lehmer1938}, pp.\ 356ff), who gave the following congruence mod $p$:

\begin{equation} \label{eq:P6a}
H_{\lfloor p/6 \rfloor} \equiv -2 \cdot q_p(2) - \frac{3}{2} \cdot q_p(3).
\end{equation}

\noindent
The fact that the vanishing of $H_{\lfloor p/6 \rfloor}$ mod $p$ is a necessary condition for the failure of the first case of FLT for the exponent $p$ is an immediate consequence of the theorems of Wieferich and Mirimanoff. In the present study, we use Lehmer's congruence in the equivalent form

\begin{equation} \label{eq:P6b}
H_{\lfloor p/6 \rfloor} \equiv -\frac{1}{2} \cdot q_p(432),
\end{equation}

\noindent
obtained from (\ref{eq:P6a}) by applying in reverse the logarithmetic and factorization rules for the Fermat quotient given by Eisenstein \cite{Eisenstein}. This expression reveals that divisors $p$ of $H_{\lfloor p/6 \rfloor}$ are instances of the vanishing of the Fermat quotient mod $p$ for composite bases, a problem which has notably been studied in the ongoing work of Richard Fischer \cite{Fischer}.

The four congruences above exhaust the cases of $H_{\lfloor p/N \rfloor}$ that can be evaluated solely in terms of Fermat quotients. Historically, the case $N$ = 5 also has its origins in this era, and was in fact introduced before that of $N$ = 6. In 1914, Vandiver \cite{Vandiver1914} proved that the vanishing of both $H_{\lfloor p/5 \rfloor}$ and $q_p(5)$ are necessary conditions for $p$ to be an exception to the first case of FLT. Unlike the four cases already considered, here the connection of the Harmonic number with FLT was discovered before any evaluation of it (beyond the definitional one) was known. It was almost eighty years later that the connection between these results would become apparent. The ingredients needed for the evaluation of this Harmonic number were presented in a 1991 paper by Williams (\cite{Williams1991}, p.\ 440), and almost simultaneously by Z.\,H.\ Sun (\cite{ZHSun1992}, pt.\ 3, Theorems 3.1 and 3.2); and though they do not write out the formula explicitly, it is clearly implied to be

\begin{equation} \label{eq:P5}
H_{\lfloor p/5 \rfloor} \equiv -\frac{5}{4} \cdot q_p(5) - \frac{5}{4} \cdot F_{p-(\frac{5}{p})}/p \pmod{p},
\end{equation}

\noindent
where $F_{p-(\frac{5}{p})}/p$ is the Fibonacci quotient (OEIS A092330), with $F$ a Fibonacci number and $\left(\frac{5}{p}\right)$ a Jacobi symbol. In light of Vandiver's theorems on $H_{\lfloor p/5 \rfloor}$ and $q_p(5)$, this result immediately established that the vanishing of the Fibonacci quotient mod $p$ was yet another criterion for the failure of the first case of FLT for the exponent $p$. The fact was announced shortly afterwards in the celebrated paper by the Sun brothers \cite{Sun+Sun1992}, which gave fresh impetus to an already vast literature on the Fibonacci quotient; and in its honor the primes $p$ which divide their Fibonacci quotient were named Wall-Sun-Sun-Primes. These remain hypothetical, as not a single instance has been found despite tests to high limits \cite{PrimeGridWallSunSunPrimeSearch}.

The remaining known formulae for the type of special Harmonic numbers in which we are interested are of much more recent origin. In some cases they were discovered simultaneously, or nearly so, by more than one researcher; and we hope we have not done injustice to any of the participants. Apart from a few published formulae which are apparently in error or underdetermined, we have the following congruences (all modulo $p$) which are undoubtably correct:

\begin{equation} \label{eq:P8}
H_{\lfloor p/8 \rfloor} \equiv -4 \cdot q_p(2) - 2 \cdot U_{p-(\frac{2}{p})}(2, -1)/p
\end{equation}

\begin{equation} \label{eq:P10}
H_{\lfloor p/10 \rfloor} \equiv -2 \cdot q_p(2) - \frac{5}{4} \cdot q_p(5) - \frac{15}{4} \cdot F_{p-(\frac{5}{p})}/p
\end{equation}

\begin{align} \label{eq:P12}
H_{\lfloor p/12 \rfloor} & \equiv -3 \cdot q_p(2) - \frac{3}{2} \cdot q_p(3) - 3 \cdot \left(\frac{3}{p}\right) \cdot U_{p-(\frac{3}{p})}(4, 1)/p \\
H_{\lfloor p/16 \rfloor} & \equiv -4 \cdot q_p(2) - 2 \cdot U_{p-(\frac{2}{p})}(2, -1)/p - 8(S-1)/p
\end{align}

\begin{equation} \label{eq:P24}
\begin{split}
H_{\lfloor p/24 \rfloor} \equiv & -4 \cdot q_p(2) -\frac{3}{2} \cdot q_p(3) - 4 \cdot U_{p-(\frac{2}{p})}(2, -1)/p \\
& - 3 \cdot \left(\frac{3}{p}\right) \cdot U_{p-(\frac{3}{p})}(4, 1)/p - 6 \cdot \left(\frac{6}{p}\right) U_{p-(\frac{6}{p})}(10, 1)/p
\end{split}
\end{equation}

\noindent
where the $\left( \frac{\cdot}{p}\right)$ are Jacobi symbols,

\begin{itemize}
\item
$F_{p-(\frac{5}{p})}/p$ is as before the Fibonacci quotient (OEIS A092330)
\item
$U_{p-(\frac{2}{p})}(2, -1)/p$ is the Pell quotient (OEIS A000129)
\item
$U_{p-(\frac{3}{p})}(4, 1)/p$ is a quotient derived from the Lucas sequence 1, 4, 15, 56, 209, \dots (OEIS A001353)
\item
$U_{p-(\frac{6}{p})}(10, 1)/p$ is a quotient derived from the Lucas sequence 1, 10, 99, 980, 9701, \dots (OEIS A004189)
\item
$\left(\frac{2}{p}\right) = (-1)^{(p^2 - 1)/8}$
\item
$\left(\frac{3}{p}\right) = (-1)^{\lfloor (p+1)/6 \rfloor}$
\item
$\left(\frac{6}{p}\right) = (-1)^{\lfloor (p+5)/12 \rfloor}$.
\end{itemize}

\noindent
and

\begin{equation} \label{eq:LerchHarmonicSplit}
S = (-1)^{\lfloor p/16 \rfloor + \lfloor p/8 \rfloor} \times
\begin{dcases}
\left( -C_{n-1} - C_{n}           + C_{n+2} \right) & \text{if $p \equiv \pm 1 \pmod{16}$} \\
\left(  C_{n-1} + C_{n} - C_{n+1}           \right) & \text{if $p \equiv \pm 3 \pmod{16}$} \\
\left(  C_{n-1}         - C_{n+1}           \right) & \text{if $p \equiv \pm 5 \pmod{16}$} \\
\left( -C_{n-1}                             \right) & \text{if $p \equiv \pm 7 \pmod{16}$},
\end{dcases}
\end{equation}

\noindent
with

\begin{subequations}
\begin{gather}
C_0 = C_1 = C_2 = C_4 = 0, C_3 = 1, C_5 = 4, C_6 = -1, C_7 = 14; \\
C_n = 8C_{n-2} - 20C_{n-4} + 16C_{n-6} - 2C_{n-8}.
\end{gather}
\end{subequations}

The result for $N$ = 8 is derived from the 1991 paper by Williams (\cite{Williams1991}, p.\ 440), with an equivalent result also appearing in Sun (\cite{ZHSun1992}, pt.\ 3, Theorem 3.3). To the best of our knowledge, at the time of the discovery of this formula in 1991 the vanishing of $H_{\lfloor p/8 \rfloor}$ modulo $p$ was not recognized as a condition for the failure of the first case of FLT for the exponent $p$, and this only became evident with the appearance of the 1995 paper of Dilcher and Skula \cite{DilcherSkula} discussed below.

The result for $N$ = 10 is due to a 1995 paper by Z.\,H.\ Sun (\cite{ZHSun1992}, pt.\ 3, Theorem 3.1). Since the vanishing of all the individual components was by then known to be a necessary condition for the failure of the first case of FLT for the exponent $p$, the same thing was immediately seen to be true for $H_{\lfloor p/10 \rfloor}$.

The result for $N$ = 12 is also derived from the 1991 paper of Williams (\cite{Williams1991}, p.\ 440), that for $N = 16$ is simplified from a 1993 paper by Z.\,H.\ Sun (\cite{ZHSun1992}, pt.\ 2, Theorem 2.1), and that for $N$ = 24 is from a 2011 paper of Kuzumaki \& Urbanowicz (\cite{KuzumakiUrbanowicz}, p.\ 139).

In their landmark joint paper of 1995, Dilcher and Skula \cite{DilcherSkula} proved among other things that the vanishing modulo $p$ of $H_{\lfloor p/N \rfloor}$ was a necessary criterion for the failure of the first case of FLT for the exponent $p$, for all $N$ from 2 to 46. This result gives retrospective interest to the evaluations of $H_{\lfloor p/8 \rfloor}$ and $H_{\lfloor p/12 \rfloor}$, and furnishes the main motivation for the present study.

All of the Fermat quotients and recurrent sequences mentioned above were evaluated in PARI using modular arithmetic.

\section{The remaining cases of $H_{\lfloor p/N \rfloor}$}

\noindent
Although the calculations in the remaining cases, where no special formula for  $H_{\lfloor p/N \rfloor}$ was known, were initially performed in the obvious way by summing of modular inverses over the required range, a better method was later found; the presentation here for the most part follows that in \cite{DobsonMatrixVariation}. The underlying idea was developed by the twin brothers Zhi-Hong Sun and Zhi-Wei Sun, who worked it out in two papers published just over a decade apart. Z.-H. Sun considered, as a special case of an even more general problem, the lacunary sum of binomial coefficients

\begin{displaymath}
T(N, m) := \sum_{\substack{j=0 \\ j \equiv 0 \pmod{N}}}^{m} \binom{m}{j}.
\end{displaymath}

\noindent
This sum is well known in the literature in its own right, having been reduced to a series in $N$ in an important paper of 1834 by C. Ramus \cite{Ramus}, who showed that

\begin{equation}
T(N, m) = \frac{1}{N} \sum_{j=0}^{N-1} \left( 2 \cos\frac{j\pi}{N} \right)^m \cos\frac{jm\pi}{N} \quad (m > 0).
\end{equation}

\noindent
This can also be written

\begin{equation}
T(N, m) = \frac{1}{N} \sum_{j=0}^{N-1} \left( 1 + \omega^j \right)^m \quad (m > 0),
\end{equation}

\noindent
where $\omega = e^{\frac{2\pi{}i}{N}}$ is a primitive $N$th root of unity. This formula was used by Hoggatt and Alexanderson \cite{HoggattAlexanderson} to derive results equivalent to the congruences (\ref{eq:P2}) through (\ref{eq:P8}), and Howard and Witt \cite{HowardWitt} extended this to the case $N=10$ (\ref{eq:P10}); in the cases $N = 5, 8, 10$ these results anticipate those cited above.

In this section, we require some supplementary notations. $H_{\lfloor p/N \rfloor}$ is the case $k=0$ of a sum studied by Lerch \cite{Lerch} and other writers,

\begin{equation} \label{eq:SkulaSum}
s(k, N) := \sum_{\substack{j=\lfloor\frac{kp}{N}\rfloor + 1\\ j \neq p}}^{\lfloor\frac{(k + 1)p}{N}\rfloor} \frac{1}{j},
\end{equation}

\noindent
where it is always assumed that $p$ is sufficiently large that $s(k, N)$ contains at least one element; the provision $j \neq p$ is necessary when $k+1=N$. It is also convenient to define an alternating version of this sum,

\begin{equation} \label{eq:SkulaSumAlternating}
s^{\ast}(k, N) := \sum_{\substack{j=\lfloor\frac{kp}{N}\rfloor + 1\\ j \neq p}}^{\lfloor\frac{(k + 1)p}{N}\rfloor} \frac{(-1)^j}{j}.
\end{equation}

Z.-H. Sun (\cite{ZHSun1992}, pt.\ 1) showed that for a prime $p$, we have

\begin{equation} \label{eq:SunFormula}
\frac{N(1 - T(N, p))}{p}  \equiv
\begin{dcases}
\medspace s(0, N) \pmod{p} & \text{if $N$ is even} \\
s^{\ast}(0, N) \equiv -s(1, 2N) \pmod{p} & \text{if $N$ is odd}.
\end{dcases}
\end{equation}

\noindent
Later, Z.-W. Sun (\cite{ZWSun2002}; see also \cite{ZWSun2008}) completed the analysis by introducing the auxiliary expression

\begin{align} \label{eq:SunDefinitionAlternating}
T^{\ast}(N, m) & := \sum_{\substack{j=0 \\ j \equiv 0 \pmod{N}}}^{m} (-1)^{j/N} \binom{m}{j} \\
               & = 2 \cdot T(2N, m) - T(N, m).
\end{align}

\noindent
Because our interest is predominantly in the harmonic numbers rather than the sums of binomial coefficients, we motivate this device in a slightly different way than in the original paper. Notice that by definition,

\begin{equation} \label{eq:SplitDefined}
s(0, N) \equiv s(0, 2N) + s(1, 2N) \pmod{p},
\end{equation}

\noindent
and by the complementarity of the elements of the residue system of a prime $p$,

\begin{displaymath}
s(k, N) \equiv -s(N - k - 1, N) \pmod{p}.
\end{displaymath}

\noindent
Let $s^\prime(k, N)$ denote the terms with odd denominators in $s(k, N)$, and $s{''}(k, N)$ the terms with even denominators. Then

\begin{align}
s^{\ast}(0, N) & \equiv s{''}(0, N) - s^\prime(k, N) \\
               & \equiv \frac{1}{2}s(0, 2N) + s{''}(N-1, N) \\
               & \equiv \frac{1}{2}s(0, 2N) + \frac{1}{2}s(N-1, 2N) \\
               & \equiv -s(1, 2N) \pmod{p},
\end{align}

\noindent
where the last line follows from Corollary 3.2 of Dilcher and Skula (\cite{DilcherSkula2011}, p. 20), which makes explicit a hint supplied in Lerch \cite{Lerch}. Thus from (\ref{eq:SplitDefined}) we have for even $N$,

\begin{equation}
s^{\ast}(0, N) \equiv s(0, 2N) - s(0, N) \pmod{p},
\end{equation}

\noindent
while for odd $N$

\begin{equation}
s(0, N) \equiv s(0, 2N) - s^{\ast}(0, N) \pmod{p}.
\end{equation}

\noindent
Substituting these last two expressions into (\ref{eq:SunFormula}), we obtain its complement,

\begin{equation} \label{eq:SunFormulaComplement}
\frac{N(1 - T^{\ast}(N, p))}{p}  \equiv
\begin{dcases}
\medspace s(0, N) \pmod{p} & \text{if $N$ is odd} \\
s^{\ast}(0, N) \equiv -s(1, 2N) \pmod{p} & \text{if $N$ is even}.
\end{dcases}
\end{equation}

\subsection{Recurrent sequence representations (I)}

\noindent
Now the sequences associated with $s(0, N)$, namely $T(N, m)$ ($N$ even) and $T^{\ast}(N, m)$ ($N$ odd), satisfy the same simple recurrence relations of order $N-1$,

\begin{equation}
a_n - \binom{N}{1}a_{n-1} + \binom{N}{2}a_{n-2} - \ldots \mp \binom{N}{N-1}a_{n-N+1} = 0
\end{equation}

\noindent
according as $N$ is even or odd; in other words,

\begin{equation}
\sum_{j=0}^{N-1} (-1)^j \binom{N}{j}a_{n-j} = 0.
\end{equation}

\noindent
Similary, the sequences associated with $s^{\ast}(0, N)$, namely $T(N, m)$ ($N$ odd) and $T^{\ast}(N, m)$ ($N$ even), satisfy recurrence relations also of order $N-1$,

\begin{equation}
a_n - \binom{N}{1}a_{n-1} + \binom{N}{2}a_{n-2} - \ldots \mp \binom{N}{N-1}a_{n-N+1} \pm 2a_{n-N} = 0
\end{equation}

\noindent
according as $N$ is even or odd; in other words,

\begin{equation}
\pm{}2a_{n-N} + \sum_{j=0}^{N-1} (-1)^j \binom{N}{j}a_{n-j} = 0.
\end{equation}

\noindent
These results are easily deduced from \cite{KonvalinaLiu}, p. 12. We will not attempt to replicate here Z.-W. Sun's proof that $T(N, m)$ can be expressed in terms of linear recurrent sequences of order not exceeding $\phi(m)/2$, where $\phi(\cdot)$ is the Euler totient function. However, we will note the advantageous decomposition

\begin{equation} \label{eq:SunFormulaSupplementEven}
s(0, 2N) \equiv s(0, N) + s^{\ast}(0, N) \pmod{p},
\end{equation}

\noindent
where by symmetry the right-hand side has the same value regardless of the parity of $N$, a value in agreement with (\ref{eq:SunFormula}).

A few examples of the characteristic polynomials in difference root form are given in Tables \ref{Table_1} and \ref{Table_2} below.

\begin{table} [H]
\begin{center}
\caption{Recurrent sequences for $T(N, m)$ for small values of $N$}
\label{Table_1}
\begin{tabular}{ r | l }
       $N$ & coefficients of the $a$'s in $a_{n} + a_{n-1} + a_{n-2} + \ldots$ = 0    \\
\hline
 2 & $1, -2$ \\
 3 & $1, -3, 3, -2$ \\
 4 & $1, -4, 6, -4$ \\
 5 & $1, -5, 10, -10, 5, -2$ \\
 6 & $1, -6, 15, -20, 15, -6$ \\
 7 & $1, -7, 21, -35, 35, -21, 7, -2$ \\
 8 & $1, -8, 28, -56, 70, -56, 28, -8$ \\
 9 & $1, -9, 36, -84, 126, -126, 84, -36, 9, -2$ \\
10 & $1, -10, 45, -120, 210, -252, 210, -120, 45, -10$ \\
\end{tabular}
\end{center}
\end{table}

\begin{table} [H]
\begin{center}
\caption{Recurrent sequences for $T^\ast(N, m)$ for small values of $N$}
\label{Table_2}
\begin{tabular}{ r | l }
       $N$ & coefficients of the $a$'s in $a_{n} + a_{n-1} + a_{n-2} + \ldots$ = 0    \\
\hline
 2 & $1, -2, 2$ \\
 3 & $1, -3, 3$ \\
 4 & $1, -4, 6, -4, 2$ \\
 5 & $1, -5, 10, -10, 5$ \\
 6 & $1, -6, 15, -20, 15, -6, 2$ \\
 7 & $1, -7, 21, -35, 35, -21, 7$ \\
 8 & $1, -8, 28, -56, 70, -56, 28, -8, 2$ \\
 9 & $1, -9, 36, -84, 126, -126, 84, -36, 9$ \\
10 & $1, -10, 45, -120, 210, -252, 210, -120, 45, -10, 2$ \\
\end{tabular}
\end{center}
\end{table}

\subsection{Matrix representations (I)}

\noindent
The most efficient representations of the sums (\ref{eq:SunFormula}) and (\ref{eq:SunFormulaComplement}) for computational purposes is believed to be in the form of matrices. While the direct evaluation of (\ref{eq:Harmonic}) as a sum of reciprocals (or modular inverses) has an algorithmetical complexity of $\frac{p \log p}{N}$, that for evaluation by matrix exponentiation is only $N^3 \log p$, and is clearly to be preferred except when $N$ is large relative to $p$, and in the present study it never is.

In what follows, whenever we speak of a term of a sequence being represented by a matrix (or the power of a matrix), we mean by the first term of the first row of the matrix (or of the resulting power); and not being aware of any standard notation for this relationship, we have used the symbol $\doteq$ to express it.

For the evaluation of recurrences of order $N-1$ there are guaranteed to be computational methods involving matrices of dimension $N-1$. However, we have preferred to use suitable matrices of dimension $N$ of a very simple type whose powers generate the terms of $T(N, m)$ and $T^{\ast}(N, m)$, feeling that the slight loss of speed with this approach was compensated for by the great ease of implementing it. The matrices for $T(\cdot, m)$ are circulant, while the matrices for $T^{\ast}(\cdot, m)$ are skew-circulant, and they differ only in the sign of the first term in the last row. We note some small examples to illustrate the pattern:

\begin{displaymath}
T(2, m) \doteq
\begin{pmatrix*}[r]
1 & 1 \\
1 & 1
\end{pmatrix*}^m
\qquad
T^{\ast}(2, m) \doteq
\begin{pmatrix*}[r]
 1 & 1 \\
-1 & 1
\end{pmatrix*}^m
\end{displaymath}

\begin{displaymath}
T(3, m) \doteq
\begin{pmatrix*}[r]
 1 & 1 & 0 \\
 0 & 1 & 1 \\
 1 & 0 & 1
\end{pmatrix*}^m
\qquad
T^{\ast}(3, m) \doteq
\begin{pmatrix*}[r]
 1 & 1 & 0 \\
 0 & 1 & 1 \\
-1 & 0 & 1
\end{pmatrix*}^m
\end{displaymath}

\begin{displaymath}
T(4, m) \doteq
\begin{pmatrix*}[r]
1 & 1 & 0 & 0 \\
0 & 1 & 1 & 0 \\
0 & 0 & 1 & 1 \\
1 & 0 & 0 & 1
\end{pmatrix*}^m
\qquad
T^{\ast}(4, m) \doteq
\begin{pmatrix*}[r]
1 & 1 & 0 & 0 \\
0 & 1 & 1 & 0 \\
0 & 0 & 1 & 1 \\
-1 & 0 & 0 & 1
\end{pmatrix*}^m
\end{displaymath}

\begin{displaymath}
T(5, m) \doteq
\begin{pmatrix*}[r]
 1 & 1 & 0 & 0 & 0 \\
 0 & 1 & 1 & 0 & 0 \\
 0 & 0 & 1 & 1 & 0 \\
 0 & 0 & 0 & 1 & 1 \\
 1 & 0 & 0 & 0 & 1
\end{pmatrix*}^m
\qquad
T^{\ast}(5, m) \doteq
\begin{pmatrix*}[r]
 1 & 1 & 0 & 0 & 0 \\
 0 & 1 & 1 & 0 & 0 \\
 0 & 0 & 1 & 1 & 0 \\
 0 & 0 & 0 & 1 & 1 \\
-1 & 0 & 0 & 0 & 1
\end{pmatrix*}^m.
\end{displaymath}

\noindent
It is much faster to calculate the powers of these matrices than to calculate the associated sequences directly, and the handling of matrix powers under modular arithmetic is very well implemented in PARI. The fact that the cost of evaluating powers of a matrix is in general proportional to the cube of its dimension $d$ is another reason for introducing (\ref{eq:SunFormulaSupplementEven}) above, because the reduction of a calculation of order $d^3$ to two calculations of order $(\frac{d}{2})^3$ involves a fourfold saving.

This method has enabled the discovery of the following solutions: for $N$ = 9, $p$ = 532199813; for $N$ = 13, $p$ = 427794751; for $N$ = 17, $p$ = 590422517.

\subsection{Recurrent sequence representations (II)}

\noindent
Several writers have observed that the recurrences for $T(N, m)$ considered above (corresponding to $s(0, N)$ for even $N$ and $s^\ast{}(0, N)$ for odd $N$) can be substantially simplified by multisecting the sequences into separate congruence classes of $m$ modulo $N$. Several OEIS entries illustrate the case $m \equiv 0 \pmod{N}$, even though they do not always explicitly state the corresponding recurrence (OEIS A007613 for $N=3$, A070775 for $N = 4$, A070782 for $N=5$, A070967 for $N=6$, A094211 for $N=7$, A070832 for $N=8$, A094213 for $N=9$, A070833 for $N=10$). We in fact avoid the case $m \equiv 0 \pmod{N}$ because it is not needed for testing prime values of $m$ and increases the order of the recurrence by 1 when $N$ is even. Table \ref{Table_3} gives the values of the coefficients of some of these sequences in difference root form. 

\begin{table} [H]
\begin{center}
\caption{Recurrent sequences for multisection of $T(N, m)$ with $m \equiv 1 \pmod{N}$ for small values of $N$}
\label{Table_3}
\begin{tabular}{ r | l }
       $N$ & coefficients of the $a$'s in $a_{n} + a_{n-1} + a_{n-2} + \ldots$ = 0    \\
\hline
 2 & $1, -4$ \\
 3 & $1, -7, -8$ \\
 4 & $1, -12, -64$ \\
 5 & $1, -21, -353, 32$ \\
 6 & $1, -38, -1691, 1728$ \\
 7 & $1, -71, -7585, 36991, 128$ \\
 8 & $1, -136, -32880, 552704, 65536$ \\
 9 & $1, -265, -139823, 6826204, 6965249, -512$ \\
10 & $1, -522, -587797, 75135226, 392963125, -3200000$ \\
\end{tabular}
\end{center}
\end{table}

\noindent
The coefficients in the second column, $-4, -7, -12, -21, -38, \ldots$ are the negatives of $2^{n-1} + n$, a shifted-index version of OEIS A005126.

Similarly, for $T^\ast(N, m)$, we avoid the case $m \equiv 0 \pmod{N}$ because it increases the order of the recurrence by 1 when $N$ is odd. Here, we must multisect the sequence according to the congruence classes of $m$ modulo $2N$ to account for the contribution of the term $2N$ in the definition of $T^\ast(N, m)$ (\ref{eq:SunDefinitionAlternating}). Table \ref{Table_4} gives the values of the coefficients of some of these sequences in difference root form.

\begin{table} [H]
\begin{center}
\caption{Recurrent sequences for multisection of $T^\ast(N, m)$ with $m \equiv 1 \pmod{2N}$ for small values of $N$}
\label{Table_4}
\begin{tabular}{ r | l }
       $N$ & coefficients in $a_{n} + a_{n-1} + a_{n-2} + \ldots$ = 0    \\
\hline
 2 & $1, 4$ \\
 3 & $1, 27$ \\
 4 & $1, 136, 16$ \\
 5 & $1, 625, 3125$ \\
 6 & $1, 2766, 172929, 64$ \\
 7 & $1, 12005, 6000099, 823543$ \\
 8 & $1, 51472, 164595808, 887132416, 256$ \\
 9 & $1, 218781, 3937153311, 360344121174, 387420489$ \\
10 & $1, 923770, 86297487035, 87413438286890, 12678013691905, 1024$ \\
\end{tabular}
\end{center}
\end{table}

% The coefficients in the second column, $.., 27, 136, 625, 2766, \ldots$ are $$.

\subsection{Matrix representations (II)}

\noindent
Here, we use the standard method of calculating terms of a recurrence by taking powers of its associated Frobenius companion matrix. The dimension $d$ of the matrix is equal to the number of terms in the tables above, minus 1 (since no column representing the first term is needed); thus $d = \lfloor (p+1)/2 \rfloor$ for $T(N, m)$, and $d = \lfloor p/2 \rfloor$ for $T^{\ast}(N, m)$. We have arbitrarily chosen to use the upper form of the companion matrix, the power of which must be multiplied by a column vector containing the $d$ starting values of each multisectioned sequence in reverse order. The required power is $\lfloor(p - 1)/N \rfloor - d + 1$ in the case of $T(N, m)$ and $\lfloor(p - 1)/(2N) \rfloor - d + 1$ in the case of $T^\ast(N, m)$, and the formula is invalid when this power is less than 0. Thus, for example, for $T^{\ast}(9, m)$, the expression takes the form

\begin{displaymath}
\begin{pmatrix*}[r]
-218781 & -3937153311 & -360344121174 & -387420489 \\
 1 & 0 & 0 & 0 \\
 0 & 1 & 0 & 0 \\
 0 & 0 & 1 & 0
\end{pmatrix*}^{\lfloor \frac{p-1}{2N} \rfloor - d + 1} \cdot v,
\end{displaymath}

\noindent
where the column vector has the value

\begin{displaymath}
v = 
\begin{dcases}
[-3399724927883844; 17199898182; -92358; 1]                   & \\
[-33347476639371456; 166522893714; -783540; 1]                & \\
[-68813711437694112; 342645201108; -1562274; 1]               & \\
[1036036074283294761; -5177777070024; 24581880; -54]          & \\
[7556243668864102341; -37885095946305; 186061536; -714]       & \\
[174286513867005385011; -876548682181710; 4443424371; -24309] & 
\end{dcases}
\end{displaymath}

\noindent
according as $p \equiv 1, 5, 7, 11, 13, 17 \bmod{18}$.

\subsection{Recurrent sequence representations (III)}

\noindent
Several writers have observed that a further reduction of the order of the recurrences for $T(N, m)$ considered above can be accomplished by instead considering $J(N, m) := N \cdot T(N, m) - 2^m$, and multisecting the resulting sequences according to the congruence classes of $m$ modulo $N$. This technique was exhibited for the case $m \equiv 0 \bmod{N}$ in a number of OEIS entries by Benoit Cloitre in 2004 (OEIS A070775 for $N = 4$, A070782 for $N=5$, A070967 for $N=6$, A070832 for $N=8$, A094213 for $N=9$). More recently, it was methodically developed in an important 2014 paper by Russell Jay Hendel \cite{Hendel}, where it may be noted that there is a misprint in equation 1.1 and in the text immediately following (see the paper's abstract for the correct version). Hendel also created the OEIS entry for these ``jump sums'' (A244608), which links to a table giving the coefficients in the recursions for $N \le 50$. It should also be mentioned that this technique is foreshadowed in the congruence for $s^\ast{}(0, 9)$ given by Z.-H. Sun in 1993 (\cite{ZHSun1992}, pt. 2, Theorem 2.7). In the examples that follow (Table \ref{Table_5}), we again avoid the case $m \equiv 0 \pmod{N}$ because it is not needed for testing prime values of $m$ and increases the order of the recurrence by 1 when $N$ is even; and with this restriction, the associated recurrences are of order $\lfloor \frac{N-1}{2} \rfloor$.

\begin{table} [H]
\begin{center}
\caption{Recurrent sequences for $J(N, m)$ for small values of $N$; cf. OEIS A244608}
\label{Table_5}
\begin{tabular}{ r | l }
       $N$ & coefficients in $a_{n} + a_{n-1} + a_{n-2} + \ldots$ = 0    \\
\hline
 3 & $1, 1$  \\
 4 & $1, 4$  \\
 5 & $1, 11, -1$  \\
 6 & $1, 26, -27$  \\
 7 & $1, 57, -289, -1$   \\
 8 & $1, 120, -2160, -256$  \\
 9 & $1, 247, -13359, -13604, 1$  \\
10 & $1, 502, -73749, -383750, 3125$  \\
11 & $1, 1013, -378283, -7682623, 1006734, 1$  \\
12 & $1, 2036, -1845522, -124221692, 126018521, 46656$  \\
13 & $1, 4083, -8689296, -1738683444, 9355414620, 107661336, -1$  \\
14 & $1, 8178, -39859401, -21957517156, 498666568799, 68871018706, -823543$  \\
\end{tabular}
\end{center}
\end{table}

\noindent
The coefficients in the second column, $0, 0, 1, 4, 11, 26, 57, 120, 247, 502, \ldots$, are the Eulerian numbers (OEIS A000295).

When $m$ is a prime $p$, the relationship between Hendel's sums $J(N, p)$ and harmonic numbers is as follows:

\begin{equation} \label{eq:HendelFormula}
\frac{N - J(N, p) - 2^p}{p} \equiv \frac{N(1 - T(N, p))}{p} \equiv
\begin{dcases}
\medspace s(0, N) \pmod{p} & \text{if $N$ is even} \\
s^{\ast}(0, N) \pmod{p} & \text{if $N$ is odd}.
\end{dcases}
\end{equation}

\noindent
There is no corresponding refinement of $T^{\ast}(N, m)$, because substituting $N \cdot T(N, m) - 2^m$ into (\ref{eq:SunDefinitionAlternating}) has no effect after the cancellation of opposite terms.

The congruence for special harmonic numbers in terms of Hendel's sums can be written

\begin{equation} \label{eq:HendelFormulaShowingFermatQuotient}
\begin{split}
\frac{N - J(N, p) - 2^p}{p} & = \frac{N - J(N, p) - 2 - (2 \cdot 2^{p-1} - 2)}{p} \\
                            & \equiv \frac{N - 2 - J(N, p))}{p} - 2q_p(2) \pmod{p}.
\end{split}
\end{equation}

\noindent
Thus, Hendel's result is equivalent to, and provides an independent derivation of, those ``classical'' congruences for $H_{\lfloor p/N \rfloor}$, with $N$ even, that contain a term $-2q_p(2)$; namely $N = 2, 6, 8, 10$ (see congruences \ref{eq:P2}, \ref{eq:P6a}, \ref{eq:P8}, \ref{eq:P10} above).

Regrettably, due to memory overruns, it proved impossible to implement Hendel's algorithm on the machinery available. However, it seems worth recording here for its possible value to other researchers.

\section{The Divisibility of Harmonic Numbers}

\noindent
It may be helpful to distinguish the purpose of the present study with that of the more general problem of the divisibility of Harmonic numbers. Eswarathasan and Levine \cite{EswarathasanLevine} note that all primes greater than 3 divide the Harmonic numbers of indices $p-1$, $p(p-1)$, and $p^2 - 1$, and Wolstenholme's theorem states that for a prime $p > 3$, $p^2 \vert H_{p-1}$. These results may be inverted to give three general rules for the divisibility of Harmonic numbers:

\begin{enumerate}
\item
$H_m$ is divisible (for even $m$) by $(m + 1)^2$ if $m + 1$ is a prime $> 3$;

\item
$H_m$ is divisible by $\frac{1 + \sqrt{4m + 1}}{2}$ if the latter expression is a prime $> 3$;

\item
$H_m$ is divisible by $\sqrt{m + 1}$ if the latter expression is a prime $> 3$.

\end{enumerate}

\noindent
Eswarathasan and Levine define \textit{harmonic} primes as primes that divide only the three aforementioned Harmonic numbers (OEIS A092101), and \textit{anharmonic} primes as those that divide additional Harmonic numbers (OEIS A092102). From the fact that we consider only cases (\ref{eq:Harmonic}) where $p \vert H_{\lfloor p/N \rfloor}$ with $N > 1$, it will be evident that any such $p$ is anharmonic, and that we are seeking a subset of anharmonic primes that divide Harmonic numbers of relatively small index; for example, from Table \ref{Table_10} we see that $N = 2$ is solved by $p = 1093$, implying that $1093 \vert H_{546}$, while $N = 46$ is solved by $p = 11731$, implying that $11731 \vert H_{255}$. In contrast, the work of Boyd \cite{Boyd} and Rogers \cite{Rogers} entails, in part, finding Harmonic numbers of large index divisible by relatively small primes, such as the case $11 \vert H_{1011849771855214912968404217247}$. The overlap between our results and theirs is slight.

Indeed, our results entail only a minority of the anharmonic primes, because they require that $p$ divide an Harmonic number $H_m$ with $m < p - 1$, and thus belong to a subset of the anharmonic primes which has sometimes been called the Harmonic irregular primes (see OEIS A092194 and the Wikipedia entry for ``Regular prime''). Although the standard definition of such a prime $p$ is that it divide an Harmonic number $H_m$ with $n < p - 1$, it may be noted that any prime satisfying this criterion must in fact divide some Harmonic number $H_m$ of smaller index $m < (p - 1)/2$, since by symmetry $H_{p-1-m} \equiv H_{m} \bmod p$. The inequality sign in this condition is strict because if $m = (p - 1)/2$, then $p$ is a Wieferich prime and by (\ref{eq:P2}) and (4) must likewise divide the smaller $H_{\lfloor p/4 \rfloor}$. But this inequality also gives a sharp upper bound on the index of the least Harmonic number divisible by an Harmonic irregular prime $p$, since among the first 999 primes we have $m = (p - 3)/2$ for $p = 29, 37, 3373$ (see OEIS A125854). In general, there does not seem to be any way of recognizing an Harmonic irregular prime other than by exhaustively testing it as a divisor of Harmonic numbers $H_m$ with $m \le (p - 3)/2$ (though one would of course not attempt to do so by actual division except for very small $p$).

In turn, our results entail only a minority of the Harmonic irregular primes, because we must have $m = \lfloor p/N \rfloor$ for some $N$. Nonetheless, for some small primes this relationship is satisfied by more than one value of $N$; for example, with $p = 137$, $p \vert H_5$, and $\lfloor p/23 \rfloor = \lfloor p/24 \rfloor = \lfloor p/25 \rfloor = \lfloor p/26 \rfloor = \lfloor p/27 \rfloor = 5$.

Taking into account the stringency of the conditions on $p$, it is perhaps unsurprising that few instances have been found where a given $p$ divides distinct Harmonic numbers of the special type under consideration. By (\ref{eq:P2}) and (4) all Wieferich primes $p$ divide both $H_{\lfloor p/4 \rfloor}$ and $H_{\lfloor p/2 \rfloor}$, but otherwise we have found only two instances where a non-Wieferich prime $p < 72,000,000$ divides two distinct Harmonic numbers of index $m$ in such a way that $m = \lfloor p/N \rfloor$ is satisfied for some $N \le 1000$ in each case: for $p = 761$, $p$ divides both $H_8$ (with $N$ = 85 through 95) and $H_{23}$ (with $N$ = 32, 33), and for $p = 845921$, $p$ divides both $H_{1011}$ (with $N$ = 836) and $H_{1524}$ (with $N$ = 555). Perhaps further solutions exist with larger $N$, but we expect such numbers to be rare.

There is, however, a sense in which our results do relate to the general problem of the divisibility of Harmonic numbers. While the formulae (\ref{eq:P2}) through (\ref{eq:P24}) are framed in terms of $p$, they may conversely be seen as divisibility conditions on $H_m$. For example, with $N=2$ (the Wieferich primes), we seek cases where the numerator of $H_m$ is divisible by a prime $2m+1$, $4m+1$, or $4m+3$. In the still unresolved case $N=5$ we ask whether it is possible for $H_m$ to be divisible by a prime of one of the forms $5m+1$, $5m+2$, $5m+3$, or $5m+4$, and in the still unresolved case $N=12$ whether by a prime of one of the forms $12m+1$, $12m+5$, $12m+7$, or $12m+11$. Incidentally, we also tested the numerators of $H_m$ for divisibility by any of these linear forms without the restriction that the divisors be prime, and found only one additional solution for $m \le 10,000$, for a divisor of the form $12n+1$: $H_{10}$ is divisible by 121, which is of course $11^2$. For the same ranges of $p$ as covered by Table \ref{Table_10} below, Table \ref{Table_6} shows, for small $N$, all linear forms for which there exists no known $p = N \cdot m + r$ ($r < N$) dividing an Harmonic number $H_m$.

\begin{table} [H]
\begin{center}
\caption{Linear forms for which no known prime $p = N \cdot m + r$ with $N \le 24$ divides an Harmonic number $H_m$; starred rows contain all possible forms of $p$; $p$ tested to at least 760,000,000}
\label{Table_6}
\begin{tabular}{ r | l }
       $N$ & $r$     \\
\hline
         2 & ---  \\
         3 & ---  \\
         4 & ---  \\
 $\star$ 5 & 1, 2, 3, 4 \\
         6 & ---          \\
         7 & 1, 3, 4, 5, 6 \\
         8 & 3, 7 \\
         9 & 1, 7, 8 \\
        10 & 1, 3 \\
        11 & 1, 2, 3, 4, 6, 7, 8, 10 \\
$\star$ 12 & 1, 5, 7, 11 \\
        13 & 1, 2, 3, 4, 5, 6, 8, 10, 11, 12 \\
        14 & 1, 3, 5, 9 \\
        15 & 1, 2, 4, 8, 11, 13, 14 \\
        16 & 5, 7, 9, 11, 13 \\
        17 & 1, 2, 3, 4, 6, 7, 8, 9, 10, 11, 12, 13, 14, 15, 16 \\
$\star$ 18 & 1, 5, 7, 11, 13, 17 \\
        19 & 1, 2, 3, 4, 5, 6, 7, 9, 10, 11, 12, 13, 14, 15, 16, 17 \\
$\star$ 20 & 1, 3, 7, 9, 11, 13, 17, 19 \\
        21 & 1, 4, 5, 8, 10, 13, 17, 19, 20 \\
        22 & 1, 3, 7, 9, 13, 15, 17, 19, 21 \\
        23 & 1, 2, 3, 4, 5, 7, 8, 9, 10, 11, 12, 13, 14, 15, 16, 17, 18, 19, 20, 21, 22 \\
        24 & 11, 13, 23 \\
\end{tabular}
\end{center}
\end{table}

\noindent
The example noted above of 121 dividing $H_{10}$ underscores the fact that if $N - 1$ is a prime, then $(N - 1)^2 \vert H_{N-2}$ and $(N - 1)^2 \equiv 1 \bmod{N}$. Thus, although there are no known \textit{prime} divisors of $H_m$ of the form $12k+1$, it is not because such divisors are \textit{algebraically} impossible.

\section{Previous Calculations}

\noindent
The following account does not aim to be complete, but we hope that it includes all significant contributions to the problem. In referring to earlier work, for the sake of brevity we focus mainly on the results of successful searches. We note the incorrectness of the statement in \cite{DilcherSkula}, pp.\ 389--390, that $H_{\lfloor p/N \rfloor} \equiv 0$ mod $p$ has a solution $p < 2,000$ for every $N$ between 2 and 46 other than 5, an error perhaps resulting from the inclusion of some vacuous sums with $p < N$. Indeed, as reported below, for $N = 18, 20, 29, 31, 43$ there are no solutions with $p < 27,580,000,000$, and for $N = 12$ there is no solution with $p < 11,545,400,000,000$.

The divisors $p$ of $H_{\lfloor p/2 \rfloor}$ and $H_{\lfloor p/4 \rfloor}$ are the Wieferich primes (OEIS A001220). The first of these, 1093, was found in 1913 by Meissner, and the second, 3511, in 1922 by Beeger. The Wieferich primes have inspired some of the most intensive numerical searches ever conducted, but no further instances have been discovered. The current search record doubtless belongs to the ongoing test by PrimeGrid \cite{PrimeGridWieferichPrimeSearch}, which had reached $p < 605,180,400,000,000,000$ as of 1 April 2016.

The divisors $p$ of $H_{\lfloor p/3 \rfloor}$ are the Mirimanoff primes (OEIS A014127). The first of these, 11, was found by Jacobi in 1828 \cite{Jacobi}, and the second, 1006003, by Kloss in 1965 \cite{Kloss}. No further instances have been discovered. The current search record appears to be $p < 970,453,984,500,000$, held by Dorais and Klyve \cite{DoraisKlyve}.

We have found no computational literature on $H_{\lfloor p/5 \rfloor}$ apart from the negative result in \cite{DilcherSkula} that it does not vanish modulo $p$ for $p < 2,000$, though inevitably, reports of negative results are more difficult to locate. One wonders why such results were not reported by Schwindt \cite{Schwindt}, since he could have obtained them as a byproduct of his study of the interval $H_{\lfloor p/5 \rfloor} - H_{\lfloor p/6 \rfloor}$. In any case, we have tested $H_{\lfloor p/5 \rfloor}$ to a much greater limit than that attainable by Schwindt's method, without finding any instance where it vanishes modulo $p$.

The divisors of $H_{\lfloor p/6 \rfloor}$ (OEIS A238201, but ignoring solutions with $p < 6$) were investigated in 1983 by Schwindt \cite{Schwindt}, whose study employed a rather na\"{\i}ve method and only reached $p < 600,000$, finding the single nontrivial case $p = 61$. We know of no further computations relating to these numbers until they were obtained in another guise as the zeros of $q_p(432)$ by Richard Fischer \cite{Fischer}, whose ongoing test, which had reached $p < 73,500,000,000,000$ as of 18 August 2016, has found two further solutions, included in Tables \ref{Table_9} and \ref{Table_10} below. The vast range studied by Fischer has not yet been completely retested, but we are continuing to replicate his calculations, as to the best of our knowledge they have been performed only once. (A similar study by Le\v{z}\'{a}k \cite{Lezak}, conducted apparently without awareness of Fischer's work, extends only to $p < 35,000,000,000$.)

\section{Computational Considerations} \label{Computational Considerations}

\noindent
The formulae involving Fermat quotients and recurrent sequences cited above are all based on exponentiation, and are essentially algorithms of the order $\log p$, while the direct evaluation of a modular inverse likewise based on exponentiation and so is of order $\log p$, but the calculation of $H_{\lfloor p/N \rfloor}$ requires the procedure to be extended over a range proportional to $p/N$. For those values of $H_{\lfloor p/N \rfloor}$ not expressible solely in terms of Harmonic numbers -- that is, apart from the cases with $N = 2, 3, 4, 6$ -- the amount of computation required for the exponentiation of the matrices is proportional to $N^3 \log p$. So while several other considerations come into play, in general terms we are (as noted above) comparing an algorithm of order $N^3 \log p$ with one of order $\frac{p \log p}{N}$. Thus, roughly speaking, as soon as we reach values of $p > N^4$, the approach involving exponentiation of matrices is to be preferred. As the largest $N$ considered in the main part of this study is 46, the break-even point for the two methods should be reached by about $p = 4,000,000$. And since the first method is proportional to $N$ while the second is inversely proportional to $N$, the greatest differences are seen in the run-times for the smallest values of $N$ for which $H_{\lfloor p/N \rfloor}$ cannot be evaluated solely in terms of Fermat quotients and require matrix exponentiation, that is, for $N = 5$. %In our implementation, by the time $p$ reached 250,000,000, the method of matrix exponentiation had become 10,000 times as fast as the method of modular inverses.

\section{The Present Calculations}

\noindent
This study does not attempt to extend existing calculations in the cases $N = 2, 3, 4, 6$, and relies on known values for those cases in the tables below. Although Cik\'{a}nek \cite{Cikanek} also gives similar conditions for the first case of FLT for $N$ up to 94 (with certain qualifications), it was decided to concentrate computing resources on the original range considered by Dilcher and Skula, with $N$ ranging from 2 to 46. The limits on $N$ in Figure 1 and in Tables \ref{Table_10} and \ref{Table_11} and are higher, as the values of $N$ between 47 and 52 happen to have fairly small $p$ as least solutions, and it was felt they may as well be recorded.

This project has been run on varying numbers of CPUs since January 2014, and it would be difficult to state with any precision the number of hours devoted to individual cases of $N$. Except for some early exploratory work, all computations were performed in PARI, using matrix exponentiation with modular arithmetic to evaluate the required recurrent sequences or matrices, as the case may be. The search for solutions of $p \vert H_{\lfloor p/N \rfloor}$ has now been carried to at least $383,950,000$ for all values of $N$ in the range, and solutions have been found for all but seven cases out of 45 (see Table \ref{Table_10} below). The search-limits so far attained vary considerably. Those cases for which special formulae exist are as follows:

\begin{table} [H]
\begin{center}
\caption{Search limits on divisors $p$ of Harmonic numbers $H_{\lfloor p/N \rfloor}$ for which a special formula is known; starred values of $N$ have no known solutions}
\label{Table_7}
\begin{tabular}{ r | r }
$N$ & limit of $p$       \\
\hline
$\star$ 5  &  13,830,000,000,000 \\
        6  &  29,010,000,000,000 \\
        8  &   6,691,500,000,000 \\
        10 &  13,830,000,000,000 \\
$\star$ 12 &  11,545,400,000,000 \\
        16 &      93,700,000,000 \\
        24 &   6,691,500,000,000 \\
\end{tabular}
\end{center}
\end{table}

\noindent
For the remaining values of $N$ the calculations were usually suspended shortly after a first solution was obtained, except for $N \le 24$, where additional data was wanted for Table \ref{Table_6}. These calculations are little different from those involving special formulae, except that the matrices are typically larger. And because the processing time for matrix exponentiation is proportional to the cube of the dimension of the matrix, the test-limits achieved for the various runs fall off as $N$ increases. The limits so far attained are as follows:

\begin{table} [H]
\begin{center}
\caption{Search limits on divisors $p$ of Harmonic numbers $H_{\lfloor p/N \rfloor}$ for which no special formula is known; starred values of $N$ have no known solutions}
\label{Table_8}
\begin{tabular}{ r | r }
$N$ & limit of $p$                                       \\
\hline
7, 14                                &    8,500,000,000   \\
9, $\star$ 18                        &  334,000,000,000   \\
11, 22                               &    2,185,000,000   \\
13, 15, 17, 19, 21, 23               &      823,100,000   \\
$\star$ 20                           &   82,880,000,000   \\
25--28, 30, 32, 33, 35--42, 44--46   &      383,950,000   \\
$\star$ 29                           &   82,879,000,000   \\
$\star$ 31                           &   68,720,000,000   \\
34                                   &      493,000,000   \\
$\star$ 43                           &   27,580,000,000   \\
\end{tabular}
\end{center}
\end{table}

\noindent
The results of the searches follow. First (Table \ref{Table_9}), we report separately those results for which the search-limits are given in Table \ref{Table_2} above, since these have particular interest and in some cases coincide with OEIS sequences. Then (Table \ref{Table_10}) we report all results combined, including those in Table \ref{Table_9}.

\section{Conclusions}

\noindent
If one considers the sequence of least divisors $p$ of Harmonic numbers $H_{\lfloor p/N \rfloor}$ for $N$ = 2 through 52 (most of which are shown in Table \ref{Table_10}), arranged in ascending order of $p$ and displayed on a logarithmetic scale (see Figure 1), it is readily seen that the points tend to describe a line of positive curvature. Thus, within the limits studied, it has been found that the sequence of least divisors $p$ grows at a superexponential rate.

\begin{figure} [ht]
	\centering
		\includegraphics{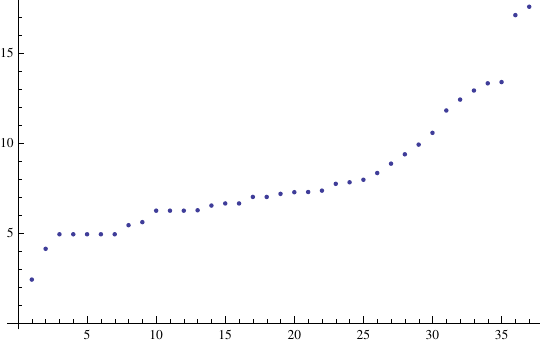}
	\caption{Least divisors $p < 27,580,000,000$ of Harmonic numbers $H_{\lfloor p/N \rfloor}$ for $N$ = 2 through 52, arranged in ascending order of $p$; horizontal axis = order of discovery; vertical axis = log $p$.}
	\label{Figure 1}
\end{figure}

Apart from the fact that the solutions for $N=2$ and $N=4$ coincide, there is no known number-theoretic reason to doubt that divisors $p$ of $H_{\lfloor p/N \rfloor}$ are randomly distributed with respect to $N$, and the data presented in Table \ref{Table_11} provide some empirical support for such a view, yielding a Pearson Correlation Coefficient between $N$ and $p$ of $-0.0166$, or nearly 0. Assuming that all cases of $N$ have at least one solution, then from a statistical point of view the scenario is a ``coupon collector's problem'' in which 44 slots (counting $N=2$ and $N=4$ as one) have to be filled with at least one solution, so that the expected number of solutions (not necessarily minimal) that must be found for all the slots to be filled is $44H_{44}$, or about 192. So far only 78 solutions (some not deriving from distinct $p$) have been found, but as the test limits for various cases of $N$ vary, we cannot infer much from this fact. However, a rough extrapolation in Mathematica\texttrademark{} from the distribution of the consecutive \textit{minimal} solutions suggests that the 44th minimal solution (if it exists) lies beyond $10^{260}$. Thus a complete solution of the problem would surely lie far beyond the present range of computability. Nonetheless, we are continuing the calculations in the hope that a few more solutions may still be found.

\section{Acknowledgements}

\noindent
We should like to thank the University of Winnipeg Library, and in particular Joffrey Abainza, for providing access to computing resources. We should also like to thank T.\,D.\ Noe for creating OEIS sequence A238201 based on our \cite{DobsonExtendedCalculations}.

\bigskip
\bigskip

\begin{table} [b]
\begin{center}
\caption{Divisors $p$ of Harmonic numbers $H_{\lfloor p/N \rfloor}$ for which a formula is known}
\label{Table_9}
\begin{tabular}{ r | l | l}
$N$ & $p$ & OEIS reference                                \\
\hline
2  & 1093, 3511                               & A001220   \\                   
3  & 11, 1006003                              & A014127   \\
4  & 1093, 3511                               & A001220   \\
5  & ---                                      & ~         \\
6  & 61, 1680023, 7308036881                  & A238201   \\
8  & 269, 8573, 1300709, 11740973, 241078561  & ~         \\
10 & 227, 17539, 4750159                      & ~         \\
12 & ---                                      & ~         \\
16 & 38723, 38993, 4292543                    & ~         \\
24 & 137, 577, 247421, 307639, 366019,        & ~         \\
~  & \quad 5262591617, 31251349243            & ~         \\
\end{tabular}
\end{center}
\end{table}

\begin{table} [tp]
\begin{center}
\caption{Divisors $p$ of Harmonic numbers $H_{\lfloor p/N \rfloor}$ for all $N$, 2 through 49}
\label{Table_10}
\begin{tabular}{ r | l || r | l }
$N$ & $p$ & $N$ & $p$ \\
\hline
2  & 1093, 3511                         & 25 & 137                 \\
3  & 11, 1006003                        & 26 & 137, 67939          \\
4  & 1093, 3511                         & 27 & 137, 23669          \\
5  & ---                                & 28 & 20101               \\     
6  & 61, 1680023, 7308036881            & 29 & ---                 \\
7  & 652913                             & 30 & 27089407            \\
8  & 269, 8573, 1300709,                & 31 & ---                 \\
~  & \quad 11740973, 241078561          & 32 & 761                 \\    
9  & 677, 6691, 532199813               & 33 & 761                 \\  
10 & 227, 17539, 4750159                & 34 & 1553                \\ 
11 & 246277, 1156457                    & 35 & 4139, 4481, 4598569 \\
12 & ---                                & 36 & 1297                \\
13 & 43214711, 427794751                & 37 & 1439, 26833         \\
14 & 2267, 6898819                      & 38 & 2473, 3527, 4047089 \\
15 & 134227                             & 39 & 407893              \\
16 & 38723, 38993, 4292543              & 40 & 509, 177553         \\
17 & 590422517                          & 41 & 509, 151883         \\
18 & ---                                & 42 & 509, 190657         \\
19 & 521, 911                           & 43 & ---                 \\
20 & ---                                & 44 & 6967, 27361         \\
21 & 1423, 5693, 5782639, 212084723     & 45 & 609221              \\
22 & 2843                               & 46 & 11731               \\
23 & 137, 264391                        & 47 & 2113                \\
24 & 137, 577, 247421, 307639, 366019,  & 48 & 2113                \\
~  & \quad 5262591617, 31251349243      & 49 & 5611043             \\
\end{tabular}
\end{center}
\end{table}

\begin{table} [hb]
\begin{center}
\caption{Least divisors $p < 27,580,000,000$ of Harmonic numbers $H_{\lfloor p/N \rfloor}$ for $N$ = 2 through 52, arranged in ascending order of $p$}
\label{Table_11}
\begin{tabular}{ r | l }
$p$ & $N$  \\
\hline
11   & 3 \\
61   & 6  \\
137  & 23, 24, 25, 26, 27 \\
227  & 10 \\
269  & 8  \\
509  & 40, 41, 42 \\
521  & 19 \\
677  & 9  \\
761  & 32, 33 \\
1093  & 2, 4 \\
1297  & 36 \\
1423  & 21 \\
1439  & 37 \\
1553  & 34 \\
2113  & 47, 48 \\
2267  & 14 \\
2473  & 38 \\
2843  & 22 \\
4127  & 51 \\
4139  & 35 \\
6967  & 44 \\
7537  & 52 \\
11731  & 46 \\
20101  & 28 \\
38723  & 16 \\
134227  & 15 \\
246277  & 11 \\
407893  & 39 \\
609221  & 45 \\
652913  & 7 \\
3337469  & 50 \\
5611043  & 49 \\
27089407  & 30 \\
43214711  & 13 \\
590422517 & 17 \\
\end{tabular}
\end{center}
\end{table}

\end{document}